\documentclass[12pt]{article}

\usepackage[top=0.75in]{geometry}
\usepackage{amsmath,amsthm,amssymb}
\usepackage{enumerate}
\usepackage{color}

\theoremstyle{definition}

\newcommand{\N}{\mathbb{N}}
\newcommand{\K}{\mathbb{K}}

\def\lf{\left\lfloor}   
\def\rf{\right\rfloor}
\def\Ker{\operatorname{Ker}}

\begin{document}

\title{Sequences: Polynomial, C-finite, Holonomic, ...}
\author{Thotsaporn Aek Thanatipanonda, Yi Zhang}

\maketitle
\thispagestyle{empty}
\vspace{-0.5in}

\begin{abstract}
Polynomial, C-finite, Holonomic are the most common
ansatz to describe the pattern of
the sequences. We propose a new ansatz called
X-recursive that generalize those we mentioned. 
We also discuss its closure properties and compare
this ansatz to another new similar ansatz from another 
paper.
\end{abstract}

The programs accompanied this article are 
\texttt{Xrec.txt} and \texttt{Guess.txt}
which can be found from the first author's website:
\texttt{thotsaporn.com}.

\section{Current Situation}
Mathematics is an art of finding patterns. 
As combinatorialists, we deal with a lot of sequences.
The process of determining the pattern of 
a given sequence is essentially important. 
The common ways to explain the sequences are by 
their recurrence relations or generating functions.
Polynomial as sequences, C-finite sequences 
(constant coefficient linear recurrences) 
and Holonomic sequences (polynomial coefficient linear recurrences) 
are the most common ansatzs as pattern searching of each sequence.
The readers are recommended to consult \cite{KP} to learn more about
these ansatzs than what we listed below.
Here we are on the expedition for a new ansatz that
can explain a bigger class of patterns.

\subsection{Polynomial as Sequences} 
The sequence $\{a(n)\}_{n=0}^{\infty}$ can be expressed as an univariate polynomial in $\K[n]$, where $\K$ is a field of characteristic of zero, and $n$ is an indeterminate, {\it i.e.}:
\[ a(n) = c_kn^k+c_{k-1}n^{k-1}+...+c_1n+c_0.\]
We call $a(n)$ a \emph{polynomial as sequence}~\cite{KP} of degree $k$.

Example 1: 
Let $a(n)=1+2+\dots+n.$ Then $a(n)=\dfrac{1}{2}n^2+\dfrac{1}{2}n.$

Since $a(n)$ is a univariate polynomial of degree $k$,
we know that the degree of $b(n) := a(n+1)-a(n)$ goes down by 1.
Applying this process for $k+1$ times, we obtain the zero sequence, {\it i.e.},
\[  (N-1)^{k+1}a(n) =0  ,  \;\ n \geq 0, \]
where the shift operator $N\cdot a(n) = a(n+1).$ 

Example 2: $a(n) := n^2-2n$ satisfies the following recurrence equation:
\[a(n+3)-3a(n+2)+3a(n+1)-a(n)=0, \;\ n \geq 0, \]
where $a(0)=0, a(1)=-1,a(2)=0.$

It follows that \textbf{the generating function} of $a(n)$ of 
polynomial as sequence of degree $k$ is
\[   \sum_{n=0}^{\infty} a(n)x^n =\dfrac{q(x)}{(1-x)^{k+1}}, \]
where $q(x)$ is some polynomial of degree at most $k.$

Some \textbf{closure properties} that are satisfied by polynomial sequences: \\
 Let $a(n)$ and $b(n)$ be a polynomial as sequence of degree $r$ and $s$. The the following properties hold,
\begin{enumerate}
\item $\{a(n) + b(n)\}_{n=0}^{\infty}$ 
is also a polynomial as sequence of degree at most $\max{(r,s)}$.
\item $\{a(n)\cdot b(n)\}_{n=0}^{\infty}$ 
is also a polynomial as sequence of degree at most $r+s.$
\item $\{\sum_{j=0}^n a(j)\}_{n=0}^{\infty}$ is also a polynomial as sequence
of degree at most $r+1.$
\end{enumerate}

\subsection{C-finite Sequences} 
The linear recurrence relation with only constant coefficients 
(\cite{KP,DZ2}), {\it i.e.}, the sequence $\{a(n)\}_{n=0}^{\infty}$ 
where there are constants  $k,c_1,\dots, c_{k-2},c_{k-1},c_{k}$ such that
\[ a(n)+c_1a(n-1)+\cdots+c_{k-1}a(n-k+1)+c_{k}a(n-k)=0 , \;\ \mbox{ for all } n \geq k.\]
In this case, we say $a(n)$ satisfies a linear recurrence of order $k.$
The polynomial sequence is a special case of C-finite sequences.

\textbf{The generating function} of C-finite sequence is
\[   \sum_{n=0}^{\infty} a(n)x^n =\dfrac{q(x)}{1+c_1x+c_2x^2+\dots+c_kx^k}, \]
where $q(x)$ is some polynomial of degree at most $k-1.$

Some examples are 
\begin{enumerate}
\item Fibonacci sequence: $a(n) = a(n-1) + a(n-2)$ ,  where $a(0)=0, a(1)=1.$
\item The sequence $\{a(n)\}_{n=0}^{\infty}$,
\[ 2, 3, 5, 9, 17, 33, 65, 129, 257, 513, 1025, ... \]
satisfies the recurrence relation
\[ a(n) = 3a(n-1)-2a(n-2),  \;\ n \geq 2.\]
and its generating function can be written as
\[  f(x) = \sum_{n=0}^{\infty} a(n)x^n = 
\dfrac{-(3x-2)}{(2x-1)(x-1)}=\dfrac{-(3x-2)}{2x^2-3x+1} .\]
\item Summation: $\displaystyle a(n)=\sum_{i=0}^{n-1}a(i), \qquad \mbox{where } a(0)=1.$ \\
This is a C-finite sequence i.e.  $a(n)=2a(n-1), \;\ n \geq 1.$
\end{enumerate}

Some \textbf{closure properties} that are satisfied by C-finite sequences: \\
 Let $a(n)$ and $b(n)$ be C-finite sequences of order $r$ and $s$. The the following properties hold,
\begin{enumerate}
\item $\{ a(n)+b(n) \}_{n=0}^{\infty}$ is a C-finite sequence of order at most $r+s,$
\item $\{ a(n)\cdot b(n) \}_{n=0}^{\infty}$ is a C-finite sequence of order at most $rs,$
\item $\{ \sum_{j=1}^{n}a(j) \}_{n=0}^{\infty}$ is a C-finite sequence of order at most $r+1,$
\item $\{ a(mn) \}_{n=0}^{\infty}, \;\ m \in \mathbb{Z^+}$ 
is a C-finite of order at most $r.$
\end{enumerate}

This closures properties are very useful for proving identities.
For example to prove that
\[  F_{2n} = 2F_nF_{n+1}-F_n^2, \;\ n \geq 0, \] 
we only define 
\[ a(n) :=  F_{2n} - 2F_nF_{n+1} + F_n^2, \;\ n \geq 0 \]
and calculate the upper bound $d$ for the order of linear recurrences of $a(n)$.
(In this case, the upper bound is 2+4+4=10.) Then 
we verify numerically that $a(i)=0, \;\ 0 \leq i \leq d$ 
to verify that $a(n)$ is the zero sequence.

\subsection{Holonomic Sequences} 
The linear recurrence relation with polynomial coefficients 
(\cite{KP}), {\it i.e.}, the sequence $\{a(n)\}_{n=0}^{\infty}$ 
where there are integer $k$ and 
polynomials  $p_0(n),p_1(n),\dots, p_{k-1}(n),p_k(n)$ 
($p_0(n) \neq 0$), each of which has degree at most $d$ such that
\[ p_0(n)a(n)+p_1(n)a(n-1)+\cdots+p_{k-1}(n)a(n-k+1)+p_k(n)a(n-k)=0 ,\] 
 for all $n \geq k.$

We say that $a(n)$ satisfies a linear recurrence of order $k.$
It is clear that the C-finite sequence is a special case of holonomic sequence.

Some examples are 
\begin{enumerate}
\item Factorial: $a(n) = n \cdot a(n-1)$ ,  where $a(0)=1.$
\item Harmonic numbers:  $H_n=1+\frac{1}{2}+\frac{1}{3}+\dots+\frac{1}{n}$
satisfies the recurrence relation,
\[ nH_{n}-(2n-1)H_{n-1}+(n-1)H_{n-2} = 0, \;\ n \geq 3   .\]
\end{enumerate}

The combinatorics community have applied this ansatz in their works for quite some time. 
But it becomes very popular since the invention of computer algebra systems like
Macsyma (early 70s) or Maple (early 80s). (Although I wish Maple had this 
command built-in.) It is the most commonly used ansatz. A lot of sequence of 
combinatorial objects fall in this class. 

\textbf{A holonomic differential equation} (aka holonomic function, D-finite function) 

Let $f(x) :=   \sum_{n=0}^{\infty} a(n)x^n$ 
be the generating function of order $k$ of
the holonomic sequence $a(n).$  
Then $f(x)$ satisfies the differential relation 
\[   q_0(x)f(x) +  q_1(x)f'(x) + \dots + q_k(x)f^{(k)}(x) =0, \]
where $q_i(x)$ are some polynomials of degree at most $d.$

We note that the order and degree of the sequence version and generating
function version do not correspond. For example, for the holonomic
sequence $a(n)$ of order $r$ and degree $d$, the corresponding
generating function has order at most $d$ and degree
at most $r+d.$

Some \textbf{closure properties} that are satisfied by holonomic sequences: \\
 Let $a(n)$ and $b(n)$ be holonomic sequences of order $r$ and $s$. 
 Then the following properties hold,
\begin{enumerate}
\item $\{ a(n)+b(n) \}_{n=0}^{\infty}$ is holonomic of order at most $r+s$,
\item  Hadamard product
$\{ a(n)\cdot b(n) \}_{n=0}^{\infty}$ is holonomic of order at most $rs$,
\item $\{ \sum_{i=0}^{n}a(i) \}_{n=0}^{\infty}$ is holonomic of order at most $r+1$,
\item Cauchy product $\{ \sum_{i=0}^n a(i)b(n-i) \}_{n=0}^{\infty}$
is holonomic of order at most $rs$,
\item $\{ a(\lf un+v\rf) \}_{n=0}^{\infty}$ is holonomic for 
any non-negative rational number $u$ and $v.$ 
\end{enumerate}

\subsection{Polynomial-Recursive Sequences} 
The sequence $\{a(n)\}_{n=0}^{\infty}$ 
where there are a polynomial with $r+1$ variables such that  
\[P(a(n), a(n-1),\dots, a(n - r)) = 0  \;\ \mbox{for all }n \geq r.\]

An example is Somos sequence, i.e. Somos-4
\[  a(n)\cdot a(n-4) - a(n-1)\cdot a(n-3) - a(n-2)^2 =0,  \;\ n \geq 5,\]
where $a(1)=a(2)=a(3)=a(4)=1.$

Michael Somos is well known for his sequences. To some surprise, 
the sequence contains only integer values. Many mathematicians 
have studied and proved the integer property of the Somos sequences.
This sequence relates directly to the theory of elliptic integral.

\section{X-Recursive Sequences}

We introduce the new ansatz which generalize the holonomic ansatz. 
We first show a motivational example.

\textbf{A hint to Somos-like sequence} 

Consider a sequence $a(n)$ generated from 
the following nonlinear relation, 
given by Michael Somos in 2014,
\[ 0 = a(n)a(n+1)a(n+3) - a(n)a(n+2)^2 - a(n+2)a(n+1)^2, \quad \mbox{for all } n \geq 0 ,  \]
where $a(0)=1, a(1)=1$ and $a(2)=2.$ 

Some of the first few terms are
\[ 1, 1, 2, 6, 30, 240, 3120, 65520, 2227680, 122522400, ... .\]
This sequence is growing too fast to be $C$-finite or holonomic, but still simple
enough for human to detect the pattern. Can you guess?

Answer:
\[ a(n) = F_{n+1} \cdot a(n-1), \qquad  a_0=1,\]
where $F_n$ is the Fibonacci sequence.
This is the sequence A003266 in OEIS website.
This example suggests the following new 
type of ansatz that might as well prove the 
integrality of Somos-like sequences. 

\subsection*{X-Recursive Sequence} 

The X-recursive sequence $\{a(n)\}_{n=0}^{\infty}$ 
 is defined to be the sequence 
where the terms satisfy a linear recurrence with 
the C-finite sequences coefficients, i.e.
\[ C_{0,n} a(n) + C_{1,n} a(n-1)+C_{2,n} a(n-2) + \cdots + C_{k,n} a(n-k) = 0 \]
where each of the sequences $C_{i,n}, \;\ 0 \leq i \leq k$ are C-finite.

We say that $a(n)$ satisfies a linear recurrence of order $k.$

We give more examples.
\begin{enumerate}
\item $ a(n) = a(n-1)+2^na(n-2), \qquad  a_0=a_1=1.$ 

\item $ a(n) = F_n \cdot a(n-1)+F_{n-1} \cdot a(n-2), \qquad  a(0)=a(1)=1.$ 

This sequence is A089126 in OEIS.

\item Summation:
$ \displaystyle a(n) = \sum_{i=1}^{n-1} F_i \cdot a(i), \qquad a_1= 1 .$

Some of the first few terms are
\[  1, 1, 2, 6, 24, 144, 1296, 18144, 399168, 13970880  .\] 
From the definition,
\[  a(n)-a(n-1) = F_{n-1}a(n-1),   \;\ n \geq 3.\] 
Therefore
\[ a(n) = (F_{n-1}+1)\cdot a(n-1)= C_n\cdot a(n-1),  \]
where $C_n = 2C_{n-1}-C_{n-3}$ with $C_3=2, C_4=3$ and $C_5=4.$
\end{enumerate}

\textbf{An X-recursive differential equation} 

Let $f(x) :=   \sum_{n=0}^{\infty} a(n)x^n$ be the generating function
of order $k$ of the X-recursive sequence $a(n)$.   
Then $f(x)$ satisfies the new relation 
\[ f(x) = \sum_{i=0}^T \left( q_{i,0}(x)f(\alpha_ix) +  q_{i,1}(x)f'(\alpha_ix) 
+ \dots + q_{i,k}(x)f^{(k)}(\alpha_ix) \right), \]
where $\alpha_i$ are the roots of $C_{i,n}$  
and $q_{i,j}(x)$ are some polynomials of degree at most $d.$

Again we note that the order and degree of the sequence version and generating
function version do not correspond directly.

We give some examples of this observation below.

Example 1: Let $a(n) = F_n \cdot a(n-1).$ Then  
\begin{align*}
f(x) &= \sum_{n=0}^{\infty} a(n)x^n 
=  \sum_{n=0}^{\infty} F_n \cdot a(n-1)x^n
= x\sum_{n=0}^{\infty} (c_1\alpha_+^n+c_2\alpha_-^n) \cdot a(n-1)x^{n-1}  \\
&= c_1'x\sum_{n=0}^{\infty} \alpha_+^n \cdot a(n)x^{n}
+c_2'x\sum_{n=0}^{\infty}\alpha_-^n \cdot a(n)x^{n} , \;\ \mbox{ (we shift index $n$ by 1)}\\
&=   c_1'xf(\alpha_+x)+ c_2'xf(\alpha_-x),
\end{align*}
where $\alpha_+$ and $\alpha_-$ are the roots of equation $x^2-x-1=0.$

For the next example, $C_n$ has a polynomial factor.

Example 2: Let $a(n) = (n+1)2^n \cdot a(n-1).$ Then 
\begin{align*}
f(x) &= \sum_{n=0}^{\infty} a(n)x^n 
=  \sum_{n=0}^{\infty} (n+1)2^n \cdot a(n-1)x^n
= 2x\sum_{n=0}^{\infty}  (n+2)2^n \cdot a(n)x^{n}  \\
&= 2x^2\sum_{n=0}^{\infty} n2^n \cdot a(n)x^{n-1}
+4x\sum_{n=0}^{\infty} 2^n \cdot a(n)x^{n} \\
&=   2x^2f'(2x)+ 4xf(2x).
\end{align*}

Next, we present \textbf{the closure properties} of X-recursive sequences. 

Let $a(n)$ and $b(n)$ be X-recursive sequences.
\begin{enumerate}
\item $\{ a(n)+b(n) \}_{n=0}^{\infty}$ is X-recursive,
\item Hadamard product
$\{ a(n)\cdot b(n) \}_{n=0}^{\infty}$ is X-recursive,
\item $\{ \sum_{i=0}^{n}a(i) \}_{n=0}^{\infty}$ is X-recursive,
\item Cauchy product $\{\sum_{i=0}^n a(i)b(n-i)\}_{n=0}^{\infty}$ 
is X-recursive.
\end{enumerate}

Properties 1-3 can be done directly by the method of undermined coefficients
They have been implemented in the program \texttt{Xrec.txt}. 
Property 4 can be proved by relating the sequence to its generating function's relationship.
The proof is the same as the holonomic sequences case mentioned in \cite[Page 142]{KP}.

All of the sequences we are talking about is in the real number field.
Things are not going so smoothly for the field with zero divisor.
This is a big deal from theoretical point of view and
we made a subsection devotes to a discussion of it. 

\subsection*{The conflict from zero 
divisors on closure properties} 

Let $\K$ be a field of characteristic zero. 
Set $\K^{\N}$ to be the ring of all sequences $\{a(n)\}_{n = 0}^\infty$ whose terms belong to $\K$. 
Assume that $N$ is the shift operator on $\K^{\N}$. 
In order to avoid sequences with only finitely many nonzero terms, 
we follow~\cite[Section 8.2]{PWZ} to take the quotient ring $\mathcal{S}(\K) := \K^{\N} / J$, where $J = \cup_{k = 0}^\infty \Ker N^k$ is the ideal 
of eventually zero sequences. 
Let $\varphi : \K^{\N} \rightarrow \mathcal{S}(\K)$ be the canonical epimorphism which maps a sequence $a \in \K^\N$ into 
its equivalence class $a + J \in \mathcal{S}(\K)$. Since $\Ker \varphi N = J$, there exist a 
unique automorphism $E$ of $\mathcal{S}(\K)$ such tat $\varphi N = E \varphi$. 
The operator $E$ is called \emph{the shift operator} on $\mathcal{S}(\K)$. 
In the following arguments, we always work in $\mathcal{S}(\K)$ and regard its elements as sequences. 
For simplicity, we use $a$ instead of the equivalence class $a + J$, and $N$ instead of $E$.  

Note that a sequence is a zero divisor in $\mathcal{S}(\K)$  if and only if it contains infinitely many zero terms and infinitely many 
nonzero terms. It implies that a sequence is a unit in $\mathcal{S}(\K)$ if and only if it is not a zero divisor. 
In order to avoid zero divisor problems, we shall work in a difference subfield~\cite{Cohn1965} of  $\mathcal{S}(\K)$ 
which contains $\K$, such as $\K(n)$. Then~\cite[Section 2]{JP} can be naturally generalized to the shift case. 
 However, to the best of my knowledge, there is no effective algorithm to test whether given finitely many sequences belong to 
a difference field or not.



 Next, we will present two examples to show that item 1 of the above closure properties only holds for some special cases, 
but not in general. 

Example 3: Assume that $a(n)$ and $b(n)$ are first-order X-recursive sequences, {\it i.e.}, 
\begin{align} \label{EQ:Xrecursive1}
C_{1, n} a(n + 1) & = C_{0, n} a(n), \\ 
D_{1, n} b(n + 1) & = D_{0, n} b(n),  \label{EQ:Xrecursive2}
\end{align} 
where $C_{0, n}, C_{1, n}, D_{0, n}, D_{1, n}$ are C-finite. 
Let us try to construct an X-recursive equation for $a(n) + b(n)$. 
Using~\eqref{EQ:Xrecursive1} and~\eqref{EQ:Xrecursive2} , we have
\begin{align}
a(n) + b(n) & = a(n) + b(n), \label{EQ:Xrecursive3} \\
C_{1, n} D_{1, n} \left(a(n + 1) + b(n + 1) \right) & = C_{0, n} D_{1, n} a(n) + C_{1, n} D_{0, n} b(n),  \label{EQ:Xrecursive4} \\
C_{1, n + 1} C_{1, n} D_{1, n+ 1} D_{1, n} \left(a(n + 2) + b(n+ 2) \right) & = C_{0, n} C_{0, n + 1} D_{1, n + 1} D_{1, n} a(n) + \nonumber \\
& \quad \ C_{1, n + 1} C_{1, n} D_{0, n + 1} D_{0, n} b(n). \label{EQ:Xrecursive5}
\end{align}
Note that we multiple $a(n + 1) + b(n + 1)$ by $C_{1, n} D_{1, n}$ because C-finite sequences are not necessarily 
units in $\mathcal{S}(\K)$. Next, let us make the following ansatz: 
\begin{multline}
x_2 C_{1, n + 1} C_{1, n} D_{1, n+ 1} D_{1, n} \left(a(n + 2) + b(n+ 2) \right) + x_1 C_{1, n} D_{1, n} \left(a(n + 1) + b(n + 1) \right) \\+ x_0 \left(a(n) + b(n) \right) = 0, \label{EQ:Xrecursive7}
\end{multline}
where $x_0, x_1, x_2$ are unknown sequences in $\mathcal{S}(\K)$.
Substituting~\eqref{EQ:Xrecursive3},~\eqref{EQ:Xrecursive4} and~\eqref{EQ:Xrecursive5} into the above equation, we get 
\begin{multline*}
\left[ C_{0, n} C_{0, n + 1} D_{1, n} D_{1, n + 1} \cdot x_2 + C_{0, n} D_{1, n} \cdot x_1 + x_0 \right] \cdot a(n) + [ C_{1, n} C_{1, n + 1} D_{0, n} D_{0, n + 1} \cdot x_2 \\ 
+ C_{1, n} D_{0, n} \cdot x_1 + x_0 ] \cdot b(n) = 0. 
\end{multline*}
Setting the coefficients of $a(n)$ and $b(n)$ to be zeros, we get the following linear equations: 
\[
 A X = 0,
\]
where
\[
A = \begin{bmatrix}
1 &  C_{0, n} D_{1, n} & C_{0, n} C_{0, n + 1} D_{1, n} D_{1, n + 1}\\
1 & C_{1, n} D_{0, n}  & C_{1, n} C_{1, n + 1} D_{0, n} D_{0, n + 1}
\end{bmatrix}, \quad \quad
X = \begin{bmatrix}
x_0 \\
x_1 \\
x_2
\end{bmatrix}.
\]
Next, we do Gaussian elimination for $A$: 
\[
A \xrightarrow{r_2 - r_1}
  \begin{bmatrix}
1 &  C_{0, n} D_{1, n} & C_{0, n} C_{0, n + 1} D_{1, n} D_{1, n + 1}\\
0 & C_{1, n} D_{0, n} - C_{0, n} D_{1, n} & C_{1, n} C_{1, n + 1} D_{0, n} D_{0, n + 1} - C_{0, n} C_{0, n + 1} D_{1, n} D_{1, n + 1}
\end{bmatrix} : = B.
\]
The second row of B corresponds to 
\begin{equation} \label{EQ:Xrecursive6}
\left( C_{1, n} D_{0, n} - C_{0, n} D_{1, n} \right) x_1 + \left( C_{1, n} C_{1, n + 1} D_{0, n} D_{0, n + 1} - C_{0, n} C_{0, n + 1} D_{1, n} D_{1, n + 1}\right) x_2  = 0. 
\end{equation}
If $C_{1, n} D_{0, n} - C_{0, n} D_{1, n} \neq 0$, then 
\[
 x_1 = \left( C_{1, n} C_{1, n + 1} D_{0, n} D_{0, n + 1} - C_{0, n} C_{0, n + 1} D_{1, n} D_{1, n + 1}\right), \quad x_2 = - \left( C_{1, n} D_{0, n} - C_{0, n} D_{1, n} \right)
\]
is a nonzero solution of~\eqref{EQ:Xrecursive6}. $C_{1, n} D_{0, n} - C_{0, n} D_{1, n} = 0$, then
\[
 x_1 = 1, \quad \quad x_2 = 0 
\]
is a nonzero solution of~\eqref{EQ:Xrecursive6}. Thus, equation~\eqref{EQ:Xrecursive6} always has a nonzero solution. From the first row of $B$, we can get the value $x_0$. 
Therefore, equation~\eqref{EQ:Xrecursive7} always has a nonzero C-finite solution. It implies that $a(n) + b(n)$ is also X-recursive. 

Example 3 shows that the sum of first-order X-recursive sequences is still X-recursive because the corresponding matrix $A$ has a special form. 
However, the next example illustrates that this might not be true for higher-order X-recursive sequences. 

Example 4: Assume that $a(n)$ and $b(n)$ are second-order X-recursive sequences, {\it i.e.}, 
\begin{align}  \label{EQ:Xrecursive8}
C_{2, n} a(n + 2) & = C_{1, n} a(n + 1) + C_{0, n} a(n), \\ 
D_{2, n} b(n + 2) & = D_{1, n} b(n + 1) + D_{0, n} b(n),  \label{EQ:Xrecursive9}
\end{align} 
where $C_{i, n}$ and $D_{i, n}$ are C-finite, $i = 0, 1, 2$. 
Let us try to construct an X-recursive equation for $a(n) + b(n)$. 
Using~\eqref{EQ:Xrecursive8} and~\eqref{EQ:Xrecursive9} , we get
\begin{align}
a(n) + b(n) & = a(n) + b(n), \label{EQ:Xrecursive10} \\
a(n + 1) + b(n + 1) & = a(n + 1) + b(n + 1),  \label{EQ:Xrecursive11} \\
\ast \cdot \left(a(n + 2) + b(n+ 2) \right) & =  \ast \cdot a(n) + \ast \cdot  a(n + 1) + \ast \cdot  b(n) + \ast \cdot  b(n + 1) \label{EQ:Xrecursive12},\\
\ast \cdot  \left(a(n + 3) + b(n+ 3) \right) & =  \ast \cdot  a(n) + \ast \cdot  a(n + 1) + \ast \cdot  b(n) + \ast \cdot  b(n + 1) \label{EQ:Xrecursive13}, \\
\ast \cdot  \left(a(n + 4) + b(n+ 4) \right) & =  \ast \cdot  a(n) + \ast \cdot  a(n + 1) + \ast \cdot  b(n) + \ast \cdot  b(n + 1) \label{EQ:Xrecursive14},
\end{align}
where $\ast$ represents a C-finite sequence. Next, let us make the following ansatz: 
\begin{multline}
x_4 \cdot  \ast \cdot  \left(a(n + 4) + b(n+ 4) \right) + x_3  \cdot  \ast \cdot  \left(a(n + 3) + b(n + 3) \right)  + x_2 \cdot  \ast \cdot  \left(a(n + 2) + b(n + 2) \right)\\
+ x_1 \left(a(n + 1) + b(n + 1)\right) + x_0 \left(a(n) + b(n) \right) = 0, \label{EQ:Xrecursive15}
\end{multline}
where $x_i$'s are unknown sequences in $\mathcal{S}(\K)$. 
Substituting equations~\eqref{EQ:Xrecursive10}-\eqref{EQ:Xrecursive14} into the above equation, we get
\begin{multline*}
 \left(\ast \cdot  x_4 + \ast \cdot  x_3 + \ast \cdot  x_2 + x_0 \right) \cdot a(n) + \left( \ast \cdot x_4 + \ast \cdot  x_3 + \ast \cdot  x_2 + x_1 \right) \cdot a(n + 1)  + \\
 \left(\ast \cdot  x_4 + \ast \cdot  x_3 + \ast \cdot  x_2 + x_0 \right) \cdot b(n) + \left( \ast \cdot x_4 + \ast \cdot  x_3 + \ast \cdot  x_2 + x_1 \right) \cdot b(n + 1)= 0.
\end{multline*}
Setting the coefficients of $a(n), a(n + 1), b(n), b(n + 1)$ to be zeros, we get the following linear equations: 
\[
 A X = 0,
\]
where
\[
A = \begin{bmatrix}
1 & 0 &  \ast & \ast & \ast \\
0 & 1 &  \ast & \ast & \ast \\
1 & 0 &  \ast & \ast & \ast \\
0 & 1 &  \ast & \ast & \ast 
\end{bmatrix}, \quad \quad
X = \begin{bmatrix}
x_0 \\
x_1 \\
x_2 \\
x_3 \\
x_4
\end{bmatrix}.
\]
Next, we do fraction-free Gaussian elimination for $A$: 
\begin{align*}
A & \xrightarrow[r_4 - r_2]{r_3 - r_1}  
\begin{bmatrix}
1 & 0 &  \ast & \ast & \ast \\
0 & 1 &  \ast & \ast & \ast \\
0 & 0 &  \ast & \ast & \ast \\
0 & 0 &  \ast & \ast & \ast 
\end{bmatrix} \\
 & \xrightarrow{\ast \cdot r_4 - \ast \cdot r_3}
 \begin{bmatrix}
1 & 0 &  \ast & \ast & \ast \\
0 & 1 &  \ast & \ast & \ast \\
0 & 0 &  \ast & \ast & \ast \\
0 & 0 &  0 & \ast & \ast 
\end{bmatrix} := B 
\end{align*}
The last row of B corresponds to 
\begin{equation} \label{EQ:Xrecursive16}
\ast \cdot x_3 + \ast \cdot x_4 = 0.
\end{equation}
Similar to the arguments in Example 3, we can always find a nonzero C-finite solution for~\eqref{EQ:Xrecursive16}. 
If we substitute this nonzero solutions into the third row of $B$, then we need to solve the following equation: 
\begin{equation} \label{EQ:Xrecursive17}
\ast_1 \cdot x_2 = \ast_2
\end{equation}
In general, $\ast_1$ might be a zero divisor in $\mathcal{S}(\K)$, such as $(-1)^n - 1$. In this case, equation~\eqref{EQ:Xrecursive17} might have no solution. 
It implies that $AX = 0$ might have no nonzero solution in $\mathcal{S}(\K)$. 
Besides, there is an algorithm~\cite{BM} to 
test whether a given C-finite sequence is a zero divisor in $\mathcal{S}(\K)$ or not. 
Even in this situation, the linear algebra approach will only be a heuristic method to test the closure properties of X-recursive 
sequence. 
\subsection*{Guessing is not easy}
The biggest problem of this ansatz is the calculation.
In practice, we must test the ansatz with the given sequences.
We demonstrate the calculation difficulties 
by showing the example where
$a(n)$ satisfies linear recurrence of order $2$ \\
\[    a(n) = C_{1,n}a(n-1) + C_{2,n}a(n-2), \;\ n \geq 2,  \]
and that $C_{1,n}$ and $C_{2,n}$ are C-finite of order 2.
Here we must solve the system of equations
\begin{align*}
a(n) &= C_{1,n}a(n-1) + C_{2,n}a(n-2), \;\ 2 \leq n \leq N, \\  
C_{1,n} &=  c_1C_{1,n-1}+c_2C_{1,n-2},   \;\ 4 \leq n \leq N, \\
C_{2,n} &=  d_1C_{2,n-1}+d_2C_{2,n-2},   \;\ 4 \leq n \leq N. 
\end{align*}

In total, there are $(N-2)+(N-4)+(N-4)=3N-10$ equations
and $2(N-2)+4=2N$ variables. Therefore we must choose
$N>10$ in order to make some sense of the guessing.
The readers may already notice that
these are also the system of non-linear equations.
It seems like the computation takes too long 
even for this simple cases.

\subsection*{Final Remarks}
At the end, we admit that we did not 
go very far with this new ansatz. 
The first problem is the slowness 
causing by solving  the system of non-linear equations 
which makes guessing very difficult.
The second problem is that the other
well known non-holonomic sequences like
Bell numbers, Somos-4,5,6 do not seem to fall 
in this class. So it is not as helpful as we might
want it to be. Nonetheless it is useful to keep
this in mind, it could be helpful when we need it one day.

\section{DD-finite Functions}

During my presentation on X-recursive sequences in 2017,
I was informed by Christoph Koutschan, my colleague from RICAM, Austria that 
the group of people at RISC leading by
Veronika Pillwein is working on a similar idea, \cite{JP}.
Instead of generalizing the sequence, they
generalized the generating function of holonomic 
functions aka D-finite functions. We will look
at this idea and compare to our 
X-recursive sequences.

\textbf{DD-finite functions} \\
Let $f(x) :=   \sum_{n=0}^{\infty} a(n)x^n.$  
$f(x)$ is called a DD-finite function of order $k$ 
if $f(x)$ satisfies the relation 
\[   q_0(x)f(x) +  q_1(x)f'(x) + \dots + q_k(x)f^{(k)}(x) =0, \]
for some D-finite functions $q_i(x)$ of degree at most a constant $d
, \;\ 0 \leq i \leq k$.

Some well known examples of these functions are

Example 1:  $e^xf(x) - f'(x) = 0$. Here $f(x) = C \cdot e^{e^x}.$
Then the coefficients $a(n)$ of $f(x)$ satisfies the relation
\[  (n+1)a(n+1) = \sum_{j=0}^n \dfrac{a(j)}{(n-j)!}   .\]
We note that if $C=e^{-1}$, $f(x)$ will be the exponential
generating function of the famous Bell numbers, $b(n)$.

By substituting $a(n) =: b(n)/n!$, we have the relation
\[  b(n+1) = \sum_{j=0}^{n} \binom{n}{j}b(j), \;\ b(0)=1.\]

Example 2: $f(x)-\sin(x)\cos(x)f'(x) = 0.$ 
This time $f(x) = C \tan(x).$
From this equation, the coefficients $a(n)$ of 
$f(x)$ satisfies the relation
\[   (1-n)a(n) = \sum_{j=1}^{n-1}jc_{n-j+1}a_j, \;\ a_0=0, a_1=1, \]
where $c_j=\dfrac{-4}{j(j-1)}c_{j-2}, \;\ c_0=0, c_1=1.$

Example 3: $((x-1)e^x+1)f(x)+x(e^x-1)f'(x) = 0.$ Here $f(x) = \dfrac{x}{e^x-1}$ 
which is the exponential generating function of the Bernoulli numbers.
From this equation, the coefficients $a(n)$ of 
$f(x)$ satisfies the relation
\[    a(n) = -\sum_{j=0}^{n-1} \dfrac{a(j)}{(n+1-j)!}.           \]

All of these three classical sequences satisfy the differential equations of order 1.

The \textbf{closure properties} of DD-finite functions 
was derived in \cite[Section 4]{JP}.
Let $f,g$ be DD-finite functions of orders $r$ and $s$. Then
\begin{enumerate}
\item $f+g$ is a DD-finite with order at most $r+s.$
\item $fg$ is a DD-finite with order at most $rs.$
\end{enumerate}

The sequence version that suggested by this ansatz is in the form
\[    \sum_{j=0}^n  c_{n,j,0}a(j) +  \sum_{j=0}^n  c_{n,j,1}a(j+1) +...
+  \sum_{j=0}^n  c_{n,j,k} a(j+k)   = 0.   \]
where, for each $l, c_{n,j,l}$ are 2-dimensional holonomic sequences in $n$ and $j.$
In practice (as of 2020), however, it is already difficult to determine whether
a giving sequence is DD-finite of order 1, i.e.
\[      a(n+1) =  \sum_{j=0}^n  c_{n,j} a(j)              .\] 
We run into the same problem
of solving system of non-linear equations as in the X-recursive case.


\end{document}